# fourier solution of two-dimensional navier stokes equation with periodic boundary conditions and incompressible flow

Logan K. Kuiper

## abstract

An approximate solution to the two dimensional Navier Stokes equation with periodic boundary conditions is obtained by representing the x and y components of fluid velocity with complex Fourier basis vectors. The chosen space of basis vectors is finite to allow for numerical calculations, but of variable size. Comparisons of the resulting approximate solutions as they vary with the size of the chosen basis vector space allow for extrapolation to an infinite basis vector space. Results suggest that such a solution, with the full infinite basis vector space and which would give the exact solution at time t, would fail for certain initial velocity configurations when initial fluid velocity and time t exceed certain finite limits.

## Introduction

The literature has many results regarding the existence of solutions for the two dimensional Navier Stokes equation. In general, the approach is to set up certain constraints on the solution for fluid velocity as a function of time t≥0 or just t=0.

Given these chosen constraints, the existence of a solution is then shown to exist, usually for all t>0 [1][2]. Wikipedia[3], citing Ladyzhenskaya, O.[4], indicates the existence of smooth and globally defined two dimensional solutions. This does not however eliminate the possibility that certain time t=0 values of velocity may disallow such a solution at some t>0. A two dimensional case with periodic boundary conditions has been presented by Zajaczkowski, W. and Zadrzynska, E.[5]

Many solutions involve taking the fourier transform of the Navier Stokes equation at the outset. This procedure is very common[6] and is very closely relate to the method being used in this work. Here, however, the approach is to simply represent the x and y components of the velocity in terms of complex fourier basis vectors. See equations [2] and [3] below.

The Navier Stokes equation is:

$$\frac{\partial v}{\partial t}+(v\cdot\Box)v=-\Box p+\alpha\Delta v \qquad (1)$$

Where v is fluid velocity, p is pressure, $\Delta$ is the Laplacian operator, $\Box$ is the gradient operator, and $\alpha$ is kinematic viscosity.

Let:

$$v_x = \sum C_{nm} e^{inx} e^{imy} \qquad (2)$$

$$v_y = \sum D_{nm} e^{inx} e^{imy} \qquad (3)$$

where n and m are both summed over the integers. Since $v_x$ and $v_y$ are real, constraints apply between the real and imaginary parts of $C_{nm}$ and also $D_{nm}$. Note that the implied periodicity of v in the x and y directions is chosen to be $2\pi$.

Incompressibility is $\Box v=0$, which along with (2) and (3) gives:

$$C_{nm} = -(m/n) D_{nm} \qquad (2a)$$

$$D_{nm} = -(n/m) C_{nm} \qquad (3a)$$

The space of integers is separated into two regions: region 1, where $|m| \geq |n|$, and region 2, where $|n| > |m|$. In region 1, m is nonzero and (3a) may be used to give Dnm from Cnm. In region 2, n is nonzero and (2a) may be used to give Cnm from Dnm.

Upon placing (2) and (3) into (1) there results:

$$BC_{nm} = -P_{nm} In \qquad (4)$$

$$BD_{nm} = -P_{nm} Im \qquad (5)$$

where

$$BC_{nm} = C_{nm} + FC_{nm}$$

$$BD_{nm} = D_{nm} + FD_{nm}$$

where

$$FC_{nm} = \sum C_{ij} C_{st} Is + \sum D_{ij} C_{st} It + \alpha C_{nm}(n^2 + m^2) \qquad (6)$$

$$FD_{nm} = \sum C_{ij} D_{st} Is + \sum D_{ij} D_{st} It + \alpha D_{nm}(n^2 + m^2) \qquad (7)$$

In these sums: i+s=n, j+t=m.

From (4) and (5) and removing $P_{nm}$:

$BC_{nm} = (n/m)BD_{nm}$              (8)

$BD_{nm} = (m/n)BC_{nm}$              (9)

(8) becomes

$(1+(n/m)^2)C_{nm} = -FC_{nm} + (n/m)FD_{nm}$     (10)

(9) becomes

$(1+(m/n)^2)D_{nm} = -FD_{nm} + (m/n)FC_{nm}$     (11)

(6), (7), and (10) become:

$(1+(n/m)^2)C_{nm} =$

$\quad -[\sum C_{ij}C_{st}Is + \sum D_{ij}C_{st}It + \alpha C_{nm}(n^2+m^2)]$

$+(n/m)[\sum C_{ij}D_{st}Is + \sum D_{ij}D_{st}It + \alpha D_{nm}(n^2+m^2)]$     (12)

Using (6), (7) and (11) instead of (10) gives:

$(1+(m/n)^2)D_{nm} =$

$\quad (m/n)[\sum C_{ij}C_{st}Is + \sum D_{ij}C_{st}It + \alpha C_{nm}(n^2+m^2)]$

$\quad +[\sum C_{ij}D_{st}Is + \sum D_{ij}D_{st}It + \alpha D_{nm}(n^2+m^2)]$     (13)

From (12):

$(1+(n/m)^2)C_{nm} = -\sum C_{ij}C_{st}Is - \alpha C_{nm}(n^2+m^2) + \ldots$

Now using

$C_{nm} = \sum_p C_{nmp} t^p$

we obtain

$(1+(n/m)^2) C_{nmp+1} = (1/(p+1))[-\sum C_{iju}C_{stv}Is - \alpha(n^2+m^2)C_{nmp} + \ldots]$  (14)

In the sums: $i+s=n$, $j+t=m$, and $u+v=p$. Using (13) instead of (12) one obtains a corresponding equation

$(1+(m/n)^2) D_{nmp+1} = \ldots\ldots$  (15)

Equations (14) and (15) are used to solve for $C_{nmp+1}$ and $D_{nmp+1}$ in regions 1 and 2 respectively. Equations (2a) and (3a) are used to get $C_{nmp+1}$ from $D_{nmp+1}$ and $D_{nmp+1}$ from $C_{nmp+1}$ respectively.

The numerical solution of equations (14) and (15) for $p=0,1,2,\ldots$ provides the solution for $v_x$ and $v_y$ as functions of $x,y$, and $t$. For $t=0$, values for $C_{nm0}$ and $D_{nm0}$ are chosen in accordance with (2a) and (3a). In the n,m plane, values for $C_{nmp}$, $p=0,N$ give rise to values for $C_{nmN+1}$, likewise for $D_{nmp}$.

One is able to observe the spread of points $C_{nmp}$ and $D_{nmp}$ in the n,m plane with increasing p, representing the spread of frequencies from initial values.

$Sum_{uv}(p) \equiv \sum_i C_{uvi} t^i$, $i=0,p$, where $p \leq 50$ is observed to determine "likely convergence"($Sum_{uv}(p)$ is constant to at least 5 significant figures for $p \geq 40$) or nonconvegence at u,v with increasing p for the determination of $C_{uv}$, this for a particular chosen t and $C_{nm0}$, $D_{nm0}$ values.

Application of the ratio test for the convergence of $Sum_{uv}(p)$ above gives the radius of convergence $R_a$ for t in $a_{uvi} t^i$ to be given by $\lim(i \to \inf) |a_{uvi}/a_{uv(i+1)}|$ and likewise for $R_b$, where $C_{uv} = a_{uv} + i b_{uv}$. This method however proves difficult to use due to the sign patterns of the terms, ++--, +-+-, etc.. The use of $Sum_{uv}$ as described above eliminates this problem.

Numerous computer runs establish "likely convergence" versus nonconvergence of $Sum_{uv}(p)$ as a function of t and the chosen initial values of $C_{nm0}$ and $D_{nmo}$. Typically $C_{nm0}$ and $D_{nmo}$ are chosen nonzero only for small |n| and |m|. Defining $v_{max}$ to be the largest of the values $|C_{nm0}|$ or $|D_{nm0}|$, "likely convergence" versus nonconvergence is established as a function of t and $v_{max}$. Results show that "likely convergence" occurs at all u,v values when it occurs. In some cases $C_{uvi} = D_{uvi} = 0$ for $i>0$, the situation with a time independent solution being chosen.

The computer program must necessarily have bounds for maximum allowed values for |n|, and |m|, |n|≤N, |m|≤N. The same bounds apply to i,j,u,s,t, and v, in (14). N was chosen to be 5,7,17, and 27(3025 points in the n,m plane). Numerous computer runs shows that "likely convergence" versus nonconvergence, mentioned above, does not usually depend

on the value for N, this for some given choice of $C_{nm0}$, $D_{nm0}$, t, and $v_{max}$. However, if "likely convergence" does occur for a particular value of N then it will also occur for any smaller value of N. If there is nonconvergence for a particular value of N, then there is also nonconvergence for any larger value for N. Thus note that if we extrapolate to a value of N larger than 27, and there is nonconvergence at N=27 for some particular choice of t and $C_{nm0}$ and $D_{nm0}$, then this would suggest nonconvergence for unbounded N.

Agreement on $Sum_{uv}(p)$ for any 2 different values of N=5,7,17,27(3025points) was to as many 6 significant figures, when "likely convergence" was obtained for both choices of N. However, for the many choices of initial $v_x$ and $v_y$ chosen, results indicate that "likely convergence" fails for sufficiently large $v_{max}$ (initial velocity) and time t. In one example, with t=5, $C_{nm0} \neq 0$ only at n,m= (-1,-1), (1,1), (2,-2), and (-2,2), "likely convergence" is obtained for all choices of N when $v_{max}$=.014 but nonconvergence for all choices of N when $v_{max}$=.14, this with $\alpha=10^{-6}$ in equ (1), $10^{-6}$ m²/sec being an approximate value for the kinematic viscosity of water. According to the previous paragraph, this result suggests that the Navier Stokes equation in two dimensions does not admit a solution in the case of $v_{max}$=.14 in the example given.

Note that the implied periodicity of all solutions for $v_x$ and $v_y$ is $2\pi$ as given in equations (2) and (3), so that the physical problem being solved is periodic in distance 2π meters when using α=$10^{-6}$.